\documentclass[]{informs4anon}

\usepackage{fancyvrb}
\usepackage{hyperref}
\usepackage[ruled]{algorithm2e}
\usepackage{pgfplotstable}
\pgfplotsset{compat=1.18}
\usepgfplotslibrary{groupplots}

\usepackage{natbib}
 \bibpunct[, ]{(}{)}{,}{a}{}{,}%
 \def\bibfont{\small}%

\TheoremsNumberedThrough
\EquationsNumberedThrough

\usepackage{xspace}
\newcommand{\MOA}[0]{\texttt{MOA.jl}\xspace}

\newcommand{\isvec}[1]{\textbf{#1}}

\definecolor{revision_color}{rgb}{0, 0, 1}
\newcommand{\revision}[1]{#1}
\newcommand{\revisionb}[1]{#1}

\begin{document}


\RUNAUTHOR{Dowson, Gandibleux, and Kof}
\RUNTITLE{MultiObjectiveAlgorithms.jl}

\TITLE{MultiObjectiveAlgorithms.jl: a Julia package for solving multi-objective optimization problems}

\ARTICLEAUTHORS{%
\AUTHOR{Oscar Dowson}
\AFF{Dowson Farms, Auckland, New Zealand, \EMAIL{oscar@dowsonfarms.co.nz}}
\AUTHOR{Xavier Gandibleux}
\AFF{Nantes Université, Nantes, France, \EMAIL{xavier.gandibleux@univ-nantes.fr}}
\AUTHOR{G{\"o}khan Kof}
\AFF{Graduate School of Sciences and Engineering, Ko\c{c} University, Istanbul, T{\"u}rkiye}
}

\ABSTRACT{
We present \texttt{MultiObjectiveAlgorithms.jl}, an open-source Julia library for solving multi-objective optimization problems written in JuMP. \texttt{MultiObjectiveAlgorithms.jl} implements \revision{ten} different solution algorithms which all rely on an iterative scalarization of the problem from a multi-objective optimization problem to a sequence of single-objective subproblems. As part of this work, we extended JuMP to support vector-valued objective functions. Because it is based on JuMP, \texttt{MultiObjectiveAlgorithms.jl} can use a wide variety of commercial and open-source solvers to solve the single-objective subproblems, and it supports problem classes \revision{including linear, integer, conic}, semi-definite, and general nonlinear.}

 

\maketitle

\section{Introduction}

This paper presents \texttt{MultiObjectiveAlgorithms.jl} (\MOA), an open-source Julia \citep{bezanson_julia_2017} library for solving multiple objective optimization models written in JuMP \citep{jump1}. Multi-objective optimization programs are well-studied in the literature; see, for example, \citet{EhrgottGandibleux2000,ehrgottMulticriteriaOptimization2005,Emmerich,EHRGOTT2025} and the references therein. We formulate a multi-objective optimization problem as:
\begin{equation}\label{eq:standard_form}
  \begin{array}{cl}
   \displaystyle \min_{\isvec{x}\in\mathbb{R}^n} & f_0(\isvec{x}) \\
   \text{s.t.} & f_i(\isvec{x}) \in S_i,\quad i = 1,\ldots,m,
  \end{array}
\end{equation}
where $f_0: \mathbb{R}^n \rightarrow \mathbb{R}^{\revision{p_0}}$ where $p_0 \ge 2$, \revision{and} $f_i: \mathbb{R}^n \rightarrow \mathbb{R}^{p_i}$ where $p_i \ge 1$, and $S_i \subset \mathbb{R}^{p_i}$ \revision{for $i = 1,\ldots,m$}. The functions $f$ and sets $S$ abstract over the various classes of mathematical optimization problems, from mixed-integer linear to semi-definite programs. For more details on this standard form, see \citet{legat2021mathoptinterface}.

\subsection{Solutions}\label{sec:intro-solutions}

\revision{Let $\mathcal{X}$ denote the set of all feasible solutions of Problem~\eqref{eq:standard_form} and let $\mathcal{Y} = \{f_0(\isvec{x}): \isvec{x} \in \mathcal{X}\}$ denote the set of corresponding objective vectors. The set $\mathcal{X}$ belongs to \textit{decision-space}, and the set $\mathcal{Y}$ belongs to \textit{objective-space}.} 

Unlike a single-objective program, where an optimal solution is reported via a single vector for $\isvec{x}$, a multi-objective program \revision{replaces} the notion of \textit{optimality} \revision{with} the notion of \textit{efficiency} \citep{ehrgottMulticriteriaOptimization2005}. \revision{Let $\le$ be defined for $\isvec{a},\isvec{b}\in\mathbb{R}^{\revision{p_0}}$ such that $\isvec{a} \le \isvec{b} \iff \isvec{b} - \isvec{a} \in \mathbb{R}^{\revision{p_0}}_+$ and $\isvec{a} \ne \isvec{b}$. Consider $\isvec{x}_1,\isvec{x}_2 \in \mathcal{X}$. If $f_0(\isvec{x}_1) \le f_0(\isvec{x}_2)$, we say $\isvec{x}_1$ \textit{dominates} $\isvec{x}_2$ and $f_0(\isvec{x}_1)$ dominates $f_0(\isvec{x}_2)$. A feasible solution $\isvec{x} \in \mathcal{X}$ is an \textit{efficient} solution if there are no feasible solutions that dominate it. The corresponding objective vector, $\isvec{y}=f_0(\isvec{x})$, is called a \textit{nondominated point}. A complete set of efficient solutions is denoted $\mathcal{X}_E \subseteq \mathcal{X}$, and the set of all nondominated points $\mathcal{Y}_N = \{f_0(\isvec{x}): \isvec{x}\in\mathcal{X}_E\}$ is called the \textit{nondominated set}. A visualization of decision- and objective-space is provided in Figure~\ref{fig:space}.}

\begin{figure}[!ht]
    \centering
    \resizebox{\textwidth}{!}{%
    \begin{tikzpicture}
        \begin{groupplot}[
            axis x line = middle, 
            axis y line = middle,
            group style = {group size = 2 by 1},
            every axis x label/.style={at={(0.5, -0.02)},anchor=north},
            every axis y label/.style={at={(-0.02, 0.5)}, anchor=south, rotate=90}
        ]
            \nextgroupplot[
                title = {Decision-space},
                xtick=\empty,
                ytick=\empty,
                xlabel = $x_1$,
                ylabel = $x_2$,
                xmin = 0,
                xmax = 5,
                ymin = 0,
                ymax = 5
            ]
                \addplot[-, fill=gray!10] coordinates{
                    (1, 2) (2, 1) (4, 2) (4, 4) (3, 4) (1, 2)
                };
                \addplot[-, mark = o, line width = 3pt] coordinates{
                    (1, 2) (2, 1) (4, 2)
                };
                \node at (axis cs:2.8,2.3) {$\mathcal{X}$};
                \node at (axis cs:1.25,1.25) {$\mathcal{X}_E$};
                 \node at (axis cs: 1, 2) [anchor = east]  {$A$};
                 \node at (axis cs: 2, 1) [anchor = north]  {$B$};
                 \node at (axis cs: 4, 2) [anchor = west]  {$C$};
            \nextgroupplot[
                title = {Objective-space},
                xtick=\empty,
                ytick=\empty,
                xlabel = $f_0^1(\isvec{x})$,
                ylabel = $f_0^2(\isvec{x})$,
                xmin = 0,
                xmax = 5,
                ymin = 0,
                ymax = 5
            ]
                \addplot[-, fill=gray!10] coordinates{
                    (1, 4) (2, 2) (4, 1) (4, 3) (3, 4) (1, 4)
                };
                \addplot[-, mark = o, line width = 3pt] coordinates{
                    (1, 4) (2, 2) (4, 1)
                };
                \node at (axis cs:3.2,2.3) {$\mathcal{Y}$};
                \node at (axis cs:1.25,2.75) {$\mathcal{Y}_N$};
                \node at (axis cs: 1, 4) [anchor = east]  {$A$};
                \node at (axis cs: 2, 2) [anchor = north east]  {$B$};
                \node at (axis cs: 4, 1) [anchor = north]  {$C$};
        \end{groupplot}
        \draw[->, dashed, line width = 1pt]
            (4,3)
            to[bend left=30]
            (12,3);
    \end{tikzpicture}}
    \caption{Visualization of decision- and objective-space for a 
    linear program with two decision variables and two objectives to minimize. The feasible set $\mathcal{X}$ in decision-space is projected to the set $\mathcal{Y}$ in objective-space via the objective function. The solid line A-B-C in decision-space is a complete set of efficient solutions $\mathcal{X}_E$. The solid line A-B-C in objective-space is the corresponding non-dominated set $\mathcal{Y}_N$.}
    \label{fig:space}
\end{figure}

Due to the existence of equivalent solutions, that is, solutions $\isvec{x}_1,\isvec{x}_2 \in \mathcal{X}_E$ such that $f_0(\isvec{x}_1) = f_0(\isvec{x}_2)$, the definition of $\mathcal{X}_E$ is refined: the \textit{maximum complete set} $\mathcal{X}_{\overline{E}}$ contains all solutions $\isvec{x}$ corresponding to all non-dominated points $\isvec{y} = f_0(\isvec{x})$ and may contain multiple $\isvec{x}$ for the same $\isvec{y}$; a \textit{complete set} $\mathcal{X}_{E}$ contains at least one $\isvec{x}$ for each non-dominated $\isvec{y}$; and a \textit{minimum complete set} $\mathcal{X}_{\underline{E}}$ contains exactly one $\isvec{x}$ for each $\isvec{y}$. Each of these sets may be infinite, for example, in the case of a multi-objective linear program in which the non-dominated set is piecewise linear. For this paper, we are concerned only with the \textit{minimum complete set} $\mathcal{X}_{\underline{E}}$. 

Some algorithms for solving Problem~\eqref{eq:standard_form} may be further restricted to finding only a \textit{representative set} of efficient solutions $\mathcal{X}_{R} \subseteq \mathcal{X}_{\underline{E}}$. \revision{A \textit{singleton representative set} is a special case when $|\mathcal{X}_R| = 1$.} \revision{Another special case concerns the notion of a \textit{supported} non-dominated point. A non-dominated point $\isvec{y}$ is \textit{supported} if $\exists \isvec{w} > 0$ for which the corresponding efficient solution $\isvec{x}$ is an optimal solution of $\min_{\isvec{x}\in\mathcal{X}}\isvec{w}^\top f_0(\isvec{x})$, otherwise it is \textit{non-supported}. Some algorithms compute} only a \textit{minimum supported set} of efficient solutions $\mathcal{X}_{\underline{SE}} \subseteq \mathcal{X}_{\underline{E}}$, which is the subset of a \textit{minimum complete set} such that all corresponding nondominated points are supported. The set of solutions \revision{computed} by an algorithm also depends on the problem type. For example, some algorithms may \revision{compute} a minimum supported set of efficient solutions for continuous linear programs, but a representative set of efficient solutions for integer or nonlinear programs. \revision{Other algorithms are applicable only for particular problem classes, for example, they may require that all decision variables are discrete.}

\subsection{Solution algorithms}
A large number of solution algorithms have been proposed for Problem~\eqref{eq:standard_form}. These algorithms can be divided into four groups along two axes according the nature of solutions returned: algorithms computing \textit{exact} or \textit{approximated} solutions, and algorithms \textit{with} or \textit{without} exploitation of preference information from the decision maker. \MOA implements \textit{exact} algorithms \textit{without} preference exploitation based on mathematical programming.

At a high level, the algorithms in \MOA solve a sequence of subproblems that are variations of the following scalar-objective problem:
\begin{equation}\label{eq:subproblem}
  \begin{array}{rl}
   \displaystyle SP(\isvec{u}, \isvec{v}) = \min_{\isvec{x}\in\mathbb{R}^n} & \isvec{u}^\top f_0(\isvec{x}) \\
   \text{s.t.} & f_i(\isvec{x}) \in S_i,\quad i = 1,\ldots,m \\
               & f_0(\isvec{x}) \le \isvec{v},
  \end{array}
\end{equation}
where $\isvec{u}\in\mathbb{R}^{\revision{p_0}}$ is a \textit{weighted-sum} vector that combines the vector-valued objective into a scalar objective, and the vector $\isvec{v}\in\mathbb{R}^{\revision{p_0}}$ places an upper bound constraint on the objective values. By iteratively changing $\isvec{u}$ and $\isvec{v}$, algorithms can find a minimum complete set. Some algorithms may add additional variables and constraints on $\isvec{x}$.

\MOA is a meta-solver. In high-level Julia code it modifies the vectors $\isvec{u}$ and $\isvec{v}$, and it uses third-party solvers to solve the resulting single-objective subproblems. Because \MOA is built on JuMP, it can use a wide variety of commercial and open-source solvers to solve the single-objective subproblems, and it supports problem classes \revision{including linear, integer, conic}, semi-definite, and general nonlinear. The solutions that \MOA returns, while subject to various theoretical guarantees, depend on the chosen algorithm and problem structure.

\revision{As an example, consider the epsilon-constraint algorithm \citep{EpsilonConstraint}. Pseudo-code for solving the bi-objective case ($p_0 = 2$) of Problem~\eqref{eq:standard_form} is provided in Algorithm~\ref{alg:epsilon}. The function $SP$ means solve Problem~\eqref{eq:subproblem} and return the optimal objective value. The function \texttt{filter\_dominated} means remove all elements of the set $\mathcal{S}$ that are dominated by another element. The algorithm is controlled by a scalar parameter $\varepsilon > 0$. When the objective functions can take only discrete integer values and $\varepsilon = 1$, the algorithm is guaranteed to find a minimum complete set, otherwise, the solution is a representative set.}

\begin{algorithm}[!ht]\label{alg:epsilon}
    $opt_1^* \leftarrow \mathrm{SP}([1, 0], [\infty, \infty])$ \\
    $opt_2^* \leftarrow \mathrm{SP}([0, 1], [\infty, \infty])$ \\
    $\mathcal{S} \leftarrow \{\}$, \ $ub_1 \leftarrow opt_1^*$, \ $ub_2 \leftarrow \infty$ \\
    \While{$ub_2 > opt_2^*$}{
         $ub_2 \leftarrow \mathrm{SP}([0, 1], [ub_1, \infty])$ and store optimal $\isvec{x}^*$\\
        $ub_1 \leftarrow ub_1 + \varepsilon$ \\
        $\mathcal{S} \leftarrow \mathcal{S} \cup \{(\isvec{x}^*, f_0(\isvec{x}^*)\}$
    }
    \Return $\texttt{filter\_dominated}(\mathcal{S})$
    \caption{Pseudo-code for the bi-objective epsilon-constraint algorithm.}
\end{algorithm}

\subsection{Contribution and outline}

The main contribution of this paper is the Julia library \texttt{MultiObjectiveAlgorithms.jl}, which we abbreviate to \MOA. A novel feature of \MOA is its modular ability to add new solution algorithms. This makes it easy for researchers to develop new solution algorithms and compare them against the existing state of the art. To enable \MOA, we extended JuMP to support vector-valued objective functions. The ease to which users can now model, solve, and analyze multi-objective optimization problems is a significant contribution of our work.

The purpose of this paper is not to provide a mathematical description of the various algorithms that are implemented in \MOA. Instead, our contribution is to make a variety of \revision{existing} algorithms available in a high quality open-source library. We refer readers \revision{interested in the algorithmic details} to the \revision{extant} literature (see Table~\ref{tab:algorithm}) and, more practically, to \revisionb{the code archive at \citet{DowsonMOA}}.

As a secondary contribution, we explain our design goals and provide a comparison to alternative software. Our hope is that our positive experiences developing \MOA provides motivation and inspiration for developers to add similar features in other algebraic modeling languages.

The rest of this paper is organized as follows. In Section~\ref{sec:example} we briefly introduce \MOA by way of a code example. \revisionb{In Section~\ref{sec:benchmarks} we provide some brief numerical benchmarks.} In Section~\ref{sec:comparison} we provide a comparison of \MOA to related work. Section~\ref{sec:principles} presents our design principles, before we conclude in Section~\ref{sec:conclusion}.

\section{Code example}\label{sec:example}

As a didactic example of the usage of \MOA, we consider a knapsack problem with two linear objectives, \revisionb{one constraint, and binary variables}:
\begin{equation}\label{eq:MOKP}
  \begin{array}{cl}
   \displaystyle\max_{\isvec{x}\in\mathbb{R}^n} & C \isvec{x} \\
   \text{s.t.} & \isvec{w}^\top \isvec{x} \le b \\
               & x_i \in \{0, 1\}, \quad i = 1,\ldots,n,
  \end{array}
\end{equation}
where $C$ is a $2$-by-$n$ matrix of costs, $\isvec{w}$ is a vector of length $n$ for the weight of each item, and $b$ is the capacity of the knapsack. The code to implement this problem in JuMP is: 

\begin{Verbatim}[fontsize=\small, frame=lines, numbers=left, xleftmargin=7mm]
using JuMP, HiGHS
import MultiObjectiveAlgorithms as MOA
function solve_knapsack(C::Matrix, w::Vector, b::Real, algorithm)
    @assert size(C) == (2, length(w))
    model = Model(() -> MOA.Optimizer(HiGHS.Optimizer))
    set_silent(model)
    set_attribute(model, MOA.Algorithm(), algorithm)
    @variable(model, x[1:length(w)], Bin)
    @objective(model, Max, C * x)
    @constraint(model, w' * x <= b)
    optimize!(model)
    assert_is_solved_and_feasible(model)
    return [value.(x; result = i) => objective_value(model; result = i)
            for i in 1:result_count(model)]
end
\end{Verbatim}
Lines 1 and 2 import the required Julia packages. Line 3 defines a Julia function with the relevant arguments (the corresponding \texttt{end} is at line 15). Line 4 performs a sanity check of the input data. Line 5 creates a new JuMP model, with \MOA as the optimizer. Because \texttt{MOA.Optimizer} is a meta-solver, it requires a mathematical programming solver as input. Here we use the HiGHS solver \citep{huangfu_2018}, which is available as \texttt{HiGHS.Optimizer} in the \texttt{HiGHS.jl} Julia package. Line 6 turns off printing. Line 7 controls which algorithm \MOA uses to solve the problem; Section~\ref{sec:supported-algorithms} provides a complete list, but examples are \texttt{MOA.EpsilonConstraint()} and \texttt{MOA.Dichotomy()}. Lines 8--12 define and solve the JuMP model using typical syntax. Note how the objective function on Line 9 is a vector. Lines 13--14 use JuMP's support for returning multiple solutions to iterate over the number of \revision{available} results and return a vector of $(\isvec{x}, \isvec{y})$ pairs corresponding to an efficient solution and its non-dominated point.

The \texttt{solve\_knapsack} function can be called as follows:

\begin{Verbatim}[fontsize=\small, frame=lines]
julia> C, w, b = rand(1:20, 2, 10), rand(1:20, 10), 50.0

julia> solve_knapsack(C, w, b, MOA.EpsilonConstraint())
4-element Vector{Pair{Vector{Float64}, Vector{Float64}}}:
 [1.0, 0.0, 1.0, 1.0, 0.0, 0.0, 1.0, 1.0, 1.0, 0.0] => [58.0, 64.0]
 [1.0, 0.0, 0.0, 1.0, 1.0, 0.0, 1.0, 1.0, 1.0, 0.0] => [65.0, 59.0]
 [1.0, 0.0, 1.0, 1.0, 1.0, 0.0, 1.0, 0.0, 1.0, 0.0] => [69.0, 56.0]
 [1.0, 0.0, 0.0, 0.0, 1.0, 1.0, 1.0, 0.0, 1.0, 0.0] => [71.0, 42.0]

julia> solve_knapsack(C, w, b, MOA.Lexicographic())
2-element Vector{Pair{Vector{Float64}, Vector{Float64}}}:
 [1.0, 0.0, 1.0, 1.0, 0.0, 0.0, 1.0, 1.0, 1.0, 0.0] => [58.0, 64.0]
 [1.0, 0.0, 0.0, 0.0, 1.0, 1.0, 1.0, 0.0, 1.0, 0.0] => [71.0, 42.0]
\end{Verbatim}

Note how the size of the solution set is algorithm-dependent. Moreover, for all algorithms, the solution returned by \MOA is a finite set of points. Thus, no algorithm can return a complete set of efficient solutions if variables are continuous. The decision to return only a finite set of points is a significant trade-off that we make for practical reasons. For more detail on this decision, see Section~\ref{sec:finite-points}.

\subsection{Supported algorithms}\label{sec:supported-algorithms}

Table~\ref{tab:algorithm} summarizes the algorithms supported by \MOA. There are currently ten algorithms implemented. Some algorithms, such as \texttt{Chalmet}, are specialized for a class of bi-objective problems. Others, such as \texttt{DominguezRios} support higher-dimensional MIP problems. The \texttt{Hierarchical} algorithm is identical to the algorithms implemented in solvers such as HiGHS \citep{huangfu_2018}. We based our implementation on the description in \citet{Hierarchical}. The \texttt{RandomWeighting} algorithm is a iterative sampling algorithm for finding a representative set. We include it in \MOA because it appears in GAMS \citep{GAMS}.

\begin{table}[!ht]
    \centering
    \resizebox{\textwidth}{!}{
    \begin{tabular}{l l c c l}
    Algorithm                  & Citation                 & \revision{Restriction} & Lines & Solution       \\
    \hline\hline
    \texttt{Chalmet}           & \citet{Chalmet}           & \revision{$p_0=2$} & \revisionb{112} & Minimum complete set            \\
    \texttt{Dichotomy}         & \citet{Dichotomy}         & \revision{$p_0=2$} & \revisionb{95} & Minimum supported set           \\
    \texttt{DominguezRios}     & \citet{DominguezRios}     & \revision{IP}      & \revisionb{220} & Minimum complete set            \\
    \texttt{EpsilonConstraint} & \citet{EpsilonConstraint} & \revision{$p_0=2$} & \revisionb{101} & Minimum complete set            \\
    \texttt{Hierarchical}      & \citet{Hierarchical}      &         &  \revisionb{92} & Singleton representative set    \\
    \texttt{KirlikSayin}       & \citet{KirlikSayin}       & \revision{IP}      & \revisionb{141} & Minimum complete set            \\
    \texttt{Lexicographic}     & \citet{Isermann1982}      &         & \revisionb{111} & Representative set              \\
    \texttt{RandomWeighting}   & \citet{GAMS}              &         &  \revisionb{61} & Representative set              \\
    \texttt{Sandwiching}       & \citet{Sandwiching}       &         & \revisionb{128} & Minimum supported set           \\
    \texttt{TambyVanderpooten} & \citet{TambyVanderpooten} & \revision{IP}      & \revisionb{231} & Minimum complete set            \\
    \end{tabular}
    }
    \caption{Summary of the algorithms supported by \MOA. \revision{The \texttt{Restriction} column indicates if the algorithm is restricted to a particular problem class: $p_0 = 2$ means the algorithm is restricted to bi-objective problems, IP means the algorithm is restricted to pure integer linear subproblems.} \textit{Lines} is the number of lines of Julia code required to implement each algorithm. The \textit{Solution} columns assumes the problem is pure linear or pure integer; for other problem classes the set of solutions may be different.}
    \label{tab:algorithm}
\end{table}

Each algorithm is relatively simple to implement, requiring \revision{less than} \revisionb{250} lines of Julia code. One reason so few lines are required is that we rely heavily on JuMP to manage the complexity of interacting with the single-objective subproblems.

Aside from the supported problem dimension, the biggest differences between the algorithms is the type of solutions that they return. Some algorithms, such as \texttt{TambyVanderpooten} can return points from a \textit{minimum complete set} for discrete problems. Other algorithms, such as \texttt{Hierarchical} return only a single point, and \texttt{Lexicographic} returns a representative set of points corresponding to all lexicographic permutations of the objective vector  (this can be disabled to return only a singleton representative set for high-dimensional problems).

\section{Benchmarks}\label{sec:benchmarks}

In this section we give a brief set of benchmarks to illustrate how \MOA can be used with different algorithms and subproblem solvers. We use the multi-objective binary knapsack instances from \citet{TambyVanderpooten}. The formulation is the same as Problem~\eqref{eq:MOKP}, except that the objective matrix $C$ is $p$-by-$n$ instead of $2$-by-$n$. We consider three problem sizes: $(p, n)=(2,200)$, $(p,n)=(3, 100)$, and $(p,n)=(4,50)$. For each problem size we use the 10 problem instances from \citet{TambyVanderpooten}. We used Julia v1.12.5 \citep{bezanson_julia_2017}, JuMP v1.30.0 \citep{jump1}, and \MOA v1.11.0. For the subproblem solvers, we use Gurobi, v13.0 \citep{gurobi} and HiGHS, v1.12 \citep{huangfu_2018}, and we use the \texttt{Dichotomy}, \texttt{EpsilonConstraint}, and \texttt{TambyVanderpooten} algorithms. When $p>2$ we study only \texttt{TambyVanderpooten}. For each case, we report the mean over the 10 instances of: the number of non-dominated solutions found; the number of single-objective subproblems solved; the total runtime in seconds; and the percentage of time spent inside the subproblem solver. All computations were carried out on a 2020 Apple M1 Mac Mini with 16 GB of RAM. We summarize the benchmark results in Table~\ref{tab:benchmark-results}. The code and data to reproduce our experiments is available at \citet{DowsonMOA}.

\begin{table}[!ht]
    \centering
    \resizebox{\textwidth}{!}{
    \revisionb{\begin{tabular}{c c l l | c c c c}
        $p$ & $n$ & Algorithm & Solver & Solutions & Subproblems & Time [s] & Time in subproblem [\%]\\
        \hline\hline
        2 & 200 & \texttt{Dichotomy} & Gurobi & 34.7 & 68.9 & 1 & 88 \\
        2 & 200 & \texttt{Dichotomy} & HiGHS & 34.9 & 68.1 & 6 & 98 \\
        2 & 200 & \texttt{EpsilonConstraint} & Gurobi & 405.4 & 813.8 & 29 & 99 \\
        2 & 200 & \texttt{EpsilonConstraint} & HiGHS & 405.4 & 813.8 & 351 & 100 \\
        2 & 200 & \texttt{TambyVanderpooten} & Gurobi & 405.4 & 818.2 & 31 & 99 \\
        2 & 200 & \texttt{TambyVanderpooten} & HiGHS & 405.4 & 818.6 & 366 & 100 \\
        3 & 100 & \texttt{TambyVanderpooten} & Gurobi & 4,090.0 & 14,724.4 & 673 & 97 \\
        4 & 50 & \texttt{TambyVanderpooten} & Gurobi & 2,240.1 & 20,570.2 & 1,042 & 68 \\
    \end{tabular}
    }}
    \caption{Benchmark results for the multiobjective knapsack problems of \citet{TambyVanderpooten} where: ``Solutions'' is the size of the non-dominated set $|Y_\mathcal{N}|$; ``Subproblems'' is the number of single-objective subproblems that were solved during the algorithm; ``Time [s]'' is the total solve time in seconds; and ``Time in subproblem [\%]'' is the percentage of the total solve time that was spent in the solve routine of the subproblems. All values are the mean across the 10 problem instances for each problem size.}
    \label{tab:benchmark-results}
\end{table}

We verified that our solutions match those reported in \citet{TambyVanderpooten}; they report the average number of non-dominated solutions as $405.4$ for $(p, n) = (2, 200)$, $4090.0$ for $(p, n) = (3, 100)$, and $2,240.1$ for $(p, n) = (4, 50)$, which matches the results for our exact algorithms. Since the \texttt{Dichotomy} algorithm finds only a minimum \textit{supported} set, it does not find all non-dominated solutions. Interestingly, our solve times are 2--5 times slower than those reported in \citet{TambyVanderpooten}, even though we solve a similar number of subproblems. Moreover, the ``Time in subproblem [\%]'' column shows that the majority of the total solve time is dominated by the time spent solving the subproblems. This means that Julia, JuMP, and \MOA are not bottlenecks in the solution process. Thus, we attribute the runtime differences to either the machine and solver used by \citet{TambyVanderpooten}, or by solver-specific options that they do not report (we used the default settings).

For a given problem size and algorithm, we observe differences in the number of subproblems solved by Gurobi and HiGHS (and also between the iteration counts reported by \citet{TambyVanderpooten}). Moreover, HiGHS finds an additional non-dominated solution in instances $9$ and $10$ of the $(p, n) = (2, 200)$ case with the \texttt{Dichotomy} algorithm. This additional solution is a valid non-dominated point, but it is not \textit{supported}. In theory, the algorithm should not find this point. However, HiGHS finds the non-supported point because the scalarized objective value is within the MIP gap of a nearby \textit{supported} non-dominated point. Thus, the different number of subproblems explored by each solver is attributable to the different sequences of iterations that can occur when there are multiple near-optimal primal solutions for a particular subproblem. 

One surprising result in Table~\ref{tab:benchmark-results} is the order-of-magnitude difference in runtime between HiGHS and Gurobi. Working with the HiGHS developers, we analyzed the instances in more detail. We verified that the performance difference was purely an artifact of the solver, and that it was not related to Julia, JuMP, \MOA, or how subproblems are efficiently cached and modified between solves. Instead, the reason for the discrepancy is simple: the subproblems are trivial mixed-integer programs with one constraint and 200 (or fewer) decision variables. HiGHS does not attempt to detect and exploit this problem structure; instead it runs generic techniques that work well on large-scale problems, but are relatively expensive on these small-scale problems.

\section{Comparison to related works}\label{sec:comparison}

Table~\ref{tab:comparison} compares \MOA to related mathematical programming software for modeling and solving multi-objective optimization problems. We compare four modeling language approaches, \MOA, AMPL \citep{AMPL}, \revision{MOCVXPY \citep{diamond2016cvxpy,mocvxpy}}, and GAMS \citep{GAMS}, and five open-source solver approaches, BenSOLVE \citep{Bensolve}, HiGHS \citep{huangfu_2018}, inner \citep{inner}, PaMILO \citep{pamilo}, and PolySCIP \citep{PolySCIP}.

Many commercial solvers (for example, COPT, CPLEX, Gurobi, Hexaly, and Xpress) have support for multiple objectives. We choose HiGHS as a single representative example of this type of solver because they all implement similar variations of the \texttt{Lexicographic} and \texttt{Hierarchical} algorithms. None support returning a set of multiple efficient solutions. To the best of our knowledge, no other popular algebraic modeling languages (for example, Pyomo \citep{bynum2021pyomo}) support multi-objective problems. MATLAB has a variety of support for multiobjective optimization, including the \texttt{paretosearch} function for returning a representative set of linear and nonlinear optimization problems \citep{MATLAB}. We omit it from Table~\ref{tab:comparison} because we could not confirm what algorithms it implements or what solvers it uses.

\begin{table}[!ht]
    \centering
    \resizebox{\textwidth}{!}{
    \begin{tabular}{r | c c c c | c c c c c}
        & \multicolumn{4}{c|}{\textit{Modeling-language-based}} & \multicolumn{5}{c}{\textit{Open-Source Solver-based}} \\
        	             & \MOA   & AMPL       & \revision{MOCVXPY} & GAMS      & BenSOLVE & HiGHS & inner & PaMILO    & PolySCIP  \\
    \hline\hline
    Programming Language & Julia  & Custom     & Python & Custom     & C        & Many  & C     & C         & C         \\
    License              & BSD-3  & Commercial & Apache & Commercial & GPL-2    & MIT   & GPL-3 & Custom    & Custom    \\
    Modeling Language    & JuMP   & AMPL       & CVXPY  & GAMS       &          &       &       &           & Zimpl     \\
    Solvers     	     & Many   & Many       & Many   & Many       & GLPK     & HiGHS & GLPK  & CPLEX/Gurobi & SCIP   \\
    Multiple Solutions   & Yes    &            & \revision{Yes} & Yes & Yes     &       & Yes   & Yes       & Yes       \\
    \hline
    \textit{Problem class} &      &            &        &            &          &       &       &           &           \\
    Linear               & Yes    & Yes        & Yes    & Yes        & Yes      & Yes   & Yes   & Yes       & Yes       \\
    Integer              & Yes    & Yes        &        & Yes        &          & Yes   &       & Yes       & Yes       \\
    Nonlinear            & Yes    & Yes        & Convex & Yes        &          &       &       & Quadratic &           \\
    \hline
    \textit{Algorithms}  &        &            &        &            &          &       &       &           &           \\
    Benson             &        &            & \revision{Yes} &    & Yes      &       & Yes   &           & Yes       \\
    B\"okler             &        &            &        &            &          &       &       & Yes       &           \\
    Chalmet              & Yes    &            &        &            &          &       &       &           &           \\
    Dichotomy            & Yes    &            &        &            &          &       &       &           &           \\
    DominguezRios        & Yes    &            &        &            &          &       &       &           &           \\
    EpsilonConstraint    & Yes    &            &        & Yes        &          &       &       &           &           \\
    Hierarchical         & Yes    & Yes        &        &            &          & Yes   &       &           &           \\
    KirlikSayin          & Yes    &            &        &            &          &       &       &           &           \\
    Lexicographic        & Yes    & Yes        &        & Yes        &          & Yes   &       &           &           \\
    RandomWeighting      & Yes    &            &        & Yes        &          &       &       &           &           \\
    Sandwiching          & Yes    &            &        & Yes        &          &       &       &           &           \\
    TambyVanderpooten    & Yes    &            &        &            &          &       &       &           &           \\
    \end{tabular}
    }
    \caption{General comparison of mathematical programming software packages for multi-objective optimization. The \textit{Algorithms} are based on the naming scheme in Table~\ref{tab:algorithm}. The exceptions are \texttt{Benson} \citep{Benson}, and \texttt{B\"okler} \citep{bokler}. The algorithms \revision{implemented by} BenSOLVE, inner, \revision{MOCVXPY,} and PolySCIP are different, but they are all variants of Benson's. PaMILO and PolySCIP have custom open-source licenses that prevent commercial use.}
    \label{tab:comparison}
\end{table}

Considering the large and well-developed body of literature that exists around multi-objective problems, mathematical programming based software packages for solving such problems are sparse and relatively feature poor. Where support does exist, most packages focus on returning a single representative solution that is some mix of a lexicographic and hierarchical blending. 

The lack of general purpose libraries means that practitioners who want to solve a multi-objective problem and return multiple solutions typically code their own implementation. In our opinion, the lack of a general purpose library for solving multi-objective optimization problems has limited the applicability of multi-objective optimization because it requires practitioners to have a reasonable level of software development proficiency. As a secondary effect, practitioners may implement only a single (and simple) solution algorithm, instead of experimenting with a range of more sophisticated solution algorithms.

The lack of general purpose libraries also hurts \textit{developers} of novel solution algorithms. In most cases, papers that develop new general purpose algorithms publish reproducible source code, but the code is often problem-specific and not intended for general purpose use. If future authors implement their algorithm in \MOA, they immediately gain the ability to model and solve problems using JuMP, and their work is also available to the wider community. We see this point as a positive, if minor, contribution of our work.

\section{Design principles}\label{sec:principles}

We based the design of \MOA on the following design principles.

\subsection{Modular algorithms}

The last decade has led to proliferation of novel solution algorithms for multi-objective optimization problems. We point to the work of \citet{DominguezRios} and \citet{TambyVanderpooten} as examples. Because algorithm design and comparison is an on-going research challenge, we wanted a design that makes writing new algorithms easy. As a secondary consideration, we wanted to make switching between algorithms simple to enable speed and solution quality comparisons between algorithms.

We achieve these goals by leveraging Julia's strength for multiple dispatch. Each algorithm subtypes the \texttt{MOA.AbstractAlgorithm} type and implements a single function, \texttt{MOA.minimize\_multiobjective}. The return value is a termination status and a list of solution points (if they exist). For example, a new algorithm could be implemented as:

\begin{Verbatim}[fontsize=\small, frame=lines]
import MultiObjectiveAlgorithms as MOA
struct NewAlgorithm <: MOA.AbstractAlgorithm end
function MOA.minimize_multiobjective(model::MOA.Optimizer, ::NewAlgorithm)
    solutions = MOA.SolutionPoint[]
    # Implement the algorithm using `model.f::MOI.AbstractVectorFunction`
    if error
        return MOI.OTHER_ERROR, nothing
    end
    return MOI.OPTIMAL, solutions
end
\end{Verbatim}
This algorithm can then be immediately used from JuMP with:
\begin{Verbatim}[fontsize=\small, frame=lines]
set_attribute(model, MOA.Algorithm(), NewAlgorithm())
\end{Verbatim}

Behind the scenes, \MOA will build the problem, convert from maximization to minimization (if necessary), and convert the returned vector of solutions into the form required by JuMP. The separation of algorithm and modeling interface simplifies writing an algorithm; each of the algorithms implemented in \MOA is less than \revisionb{250} lines of code.

\subsection{A solution is a finite set of points}\label{sec:finite-points}

Unlike the single-objective optimization literature, the multi-objective algorithm literature has not converged on a single definition of what constitutes a ``solution.'' 
As described in Section~\ref{sec:intro-solutions}, there are at least three common definitions: the \textit{maximum complete set}, a \textit{complete set}, and a \textit{minimum complete set}. These sets may be finite and discrete or continuous and infinite.

In \MOA, we choose to ignore this complexity. We do not require that an algorithm returns a minimum complete set. We define a solution as a finite set of vectors $\isvec{x}$ with corresponding objective vectors $\isvec{y} = f_0(\isvec{x})$. \revision{We order solutions lexicographically based on the objective vector and the objective sense; the first solution is best.} We make no quality assertions about the returned list, other than that each objective vector is non-dominated. We make this decision for pragmatic reasons: it makes the library simpler to implement, and it provides most users with most of what they want in practice, which is the ability to find a representative set of diverse solutions to their problem. The trade-off is that it requires users to have some theoretical knowledge to choose an appropriate solution algorithm for their problem type, and to understand whether the returned solution is a representative set or a minimum complete set.

\subsection{A single objective sense and a vector-valued objective}

When designing \MOA, we had robust discussions about two alternate approaches: is a multi-objective optimization problem defined by a set of independent scalar objectives, each of which may have a different objective sense, or is it a vector-valued function with a single objective sense? Many users would prefer the former, since it allows them to, for example, minimize one objective and maximize another. The latter requires them to normalize all scalar objectives into the same objective sense by multiplying maximization objectives by $-1$.

Because \MOA is based on JuMP, which itself relies on the MathOptInterface standard form \citep{legat2021mathoptinterface}, we chose to implement the library as a single objective sense and a vector-valued objective function. Although more restrictive from a user-perspective, this significantly simplified the implementation of \MOA, and we have not found this decision to be a burden in practice. Note that our decision is purely an artifact of our decision to build on JuMP; GAMS \revision{and MOCVXPY}, for example, allow \revision{a} different objective sense \revision{for each scalar objective}.

\section{Conclusions}\label{sec:conclusion}

This paper has presented \texttt{MultiObjectiveAlgorithms.jl} (\MOA), an open-source Julia package for solving multi-objective optimization problems written in JuMP. \revisionb{Readers are directed to \citet{DowsonMOA} for the archived source code and to \url{https://github.com/jump-dev/MultiObjectiveAlgorithms.jl} for the latest version.}

\MOA implements ten different solution algorithms, which all rely on an iterative scalarization of the problem from a multi-objective optimization problem to a sequence of single-objective subproblems. Each algorithm is relatively simple to implement, requiring \revision{less than} \revisionb{250} lines \revision{of code}, and the modular design of the library makes it easy to add new algorithms. Because it is based on JuMP, \MOA can use a wide variety of commercial and open-source solvers to solve the single-objective subproblems, and it supports problem classes ranging from linear, to conic, semi-definite, and general nonlinear. 

The only gap between recent advances in algorithms for multi-objective optimization and the current design of the library is our decision to return a finite set of points as the solution. As future work, we would like to extend \MOA to support solutions such as nondominated facets and the associated open and closed points. The correct implementation of this is an open question, particularly with how to represent these solutions in JuMP.

To conclude, we hope that our positive experience developing \MOA serves as motivation and design inspiration for others to add similar features to algebraic modeling languages that do not currently support multi-objective optimization problems.

\ACKNOWLEDGMENT{
\MOA is the successor to the \texttt{vOptSolver} project \citep{vOptSolver2017}. The development of vOptSolver started in the ANR/DFG-14-CE35-0034-01 research project vOpt (2015--2019).
}

\bibliographystyle{informs2014}
\bibliography{main}

\end{document}